\documentclass[11pt,a4paper]{article}
\usepackage{fullpage}
\usepackage{graphicx}
\usepackage{subfigure}
\usepackage{amsmath, amssymb}
\date{}%leave empty

\newcommand{\R}{\mathbb{R}}

\def\abs#1{\mid#1\mid}

\def\div{{\rm div}}

\def\1{1\hspace{-.6ex} {\rm{I}}}

\newtheorem{rem}{Remark}
\newtheorem{thm}{Theorem}

\title{A kinetic scheme for unsteady pressurised flows in closed water pipes}

\author{C. Bourdarias$^{1}$ and S. Gerbi$^{1}$ \\[1mm]
$^{1}$Universit\'{e} de Savoie,
Laboratoire de Math\'{e}matiques\\
73376 Le Bourget-du-Lac Cedex, France.\\[1mm]
{\small e-mails: Christian.Bourdarias@univ-savoie.fr, Stephane.Gerbi@univ-savoie.fr}
}
\begin{document}
\maketitle
\begin{abstract}
The aim of this paper is to present a kinetic numerical scheme for the computations of transient
pressurised flows in closed water
pipes. Firstly, we detail the mathematical model written as a conservative hyperbolic partial
differentiel system of equations,
and then we recall how to obtain the corresponding kinetic formulation. Then we build the
kinetic scheme ensuring an
upwinding of the source term due to the topography  performed in a close manner described by
Perthame et al.
\cite{PS01, BPV03} using an energetic balance at microscopic level.
The validation is lastly performed in the case of a water hammer in an uniform pipe:
we compare the numerical results provided by an industrial code used at EDF-CIH (France), which
solves the Allievi equation
(the commonly used equation for pressurised flows in pipes) by the method of characteristics, with
those of the kinetic scheme. It appears that they are in a very good agreement.
\end{abstract}%%%%%%%%%%%%%%%%%%%%%%%%%%%%%%%%%%%%%%%%%%%%%%%%%%%%%%%%%%%%%%%%%%%%%%%%%%%%%%
\section{Introduction}
%%%%%%%%%%%%%%%%%%%%%%%%%%%%%%%%%%%%%%%%%%%%%%%%%%%%%%%%%%%%%%%%%%%%%%%%%%%%%%
The work presented in this article is the first step in a more general project: the mode\-li\-sation
of unsteady mixed
water flows in open channels and in pipes, its kinetic formulation and its numerical resolution by a
kinetic scheme.

Since we are interested in flows occuring in closed pipes, it may happen that some parts of the flow
are free-surface (this means that only a part of the cross-section of the pipe is filled)
and other parts are pressurised  (this means that all the cross-section of the pipe is filled).
Let us thus recall the current and previous works about mixed flows in closed water pipes
by a partial state of the art review.
%----------------------------------------------------------------------------
Cunge and Wegner \cite{CW64}  studied the pressurised flow in a pipe as if
it were a free-surface flow by assuming a narrow slot to exist in
the upper part of the tunnel, the width of the slot being calculated
to provide the correct sonic speed. This approach has been credited to Preissmann.
Later, Cunge \cite{C66} conducted a study of translation waves in a power canal containing a
series of transitions, including a siphon. Pseudoviscosity methods were employed to describe
the movement of bores in open-channel reaches.
Wiggert \cite{W72} studied the transient flow phenomena and his analytical considerations
included open-channel surge equations that were solved by the method of characteristics.
He subjected it to subcritical flow conditions. His solution resulted
from applying a similarity between the movement of a hydraulic bore and an interface (that is,
a surge front wave). Following Wiggert's model, Song, Cardle and Leung \cite{SCL83}
developed two mathematical models of unsteady free-surface/pressurised flows
using the method of characteristics (specified time and space) to compute flow
conditions in two flow zones. They showed that the pressurised phenomenon is a
dynamic shock requiring a full dynamic treatment even if inflows and other boundary conditions
change very slowly. However the Song models do not include the bore presence in the free-surface zone.
Hamam and McCorquodale  \cite{HM82} proposed a rigid water column approach to model the mixed flow pressure transients.
This model assumes a hypothetical stationary bubble across compression and expansion processes.
Li and McCorquodale \cite{LM99} extended the rigid water column approach to allow for the transport
of the trapped air bubble. Recently Fuamba \cite{M02} proposed a model for the transition from a free surface flow to a
pressurised. He wrote the conservation of mass, momentum and energy through the transition point
and proposed a laboratory validation of his model.
In the last few years, numerical models mainly based on the Preissmann
slot technique have been developed to handle the flow transition in sewer systems.
Implementing the Preissmann slot technique has the advantage of
using only one flow type (free-surface flow) throughout the whole
pipe and of being able to easily quantify the pressure head when pipes pressurise.
Let us specially mention the work of Garcia-Navarro, Alcrudo and Priestley \cite{GAP94}
in which an implicit method based on the characteristics has been proposed.

The Saint Venant equations, which  are  written in a conservative form, are usually used
to describe free surface flows of water in  open channels. As said before, they are also used in the context
of mixed flows (i.e. either free surface or pressurized) using the artifice of the Preissmann
slot \cite{SWB98},\cite{CSZ97}. On the other hand, the commonly used model to describe pressurized flows in
pipes is the system of the Allievi equations \cite{SWB98}. This system
of 1st order partial differential equations cannot be written under a conservative form
since this model is derived by neglecting some acceleration terms. This non conservative formulation is
not appropriate for a kinetic interpretation of the transition between
the two types of flows since we are not able to write conservations of appropriate quantities such as momentum and energy.

Then, it appears that a conservative model and a kinetic interpretation of it which describes pressurised flows in closed
water pipes could be of great interest.

The model used in this article to describe pressurised flows in closed water pipe is very closed to
the Shallow Water equations, and has been established by the authors in \cite{BG08}. A second order
well-balanced finite volume scheme was therein presented.
We will recall in section \ref{model} the main features of this previous work.

Another approach for the numerical resolution of Shallow Water equations is
to use a kinetic formulation \cite{PS01, BPV03}. The corresponding scheme appears to have
interesting theoretical properties:
the scheme preserves the still water steady state and involves a conservative in-cell entropy
inequality.
Moreover, this type of numerical approximation leads to an easy implementation. The present
mo\-delisation of pressurised flows is formally very close to the Shallow Water equations and it may be very interesting to
propose a kinetic formulation and thus to construct a kinetic scheme.

The model for the unsteady mixed water flows in closed water pipes and a finite volume
discretisation has been previoulsy studied by the authors \cite{BG07} and a kinetic formulation  has
been proposed in \cite{BGG08}.
We will recall in section \ref{model} the main results and the properties of this kinetic
formulation that will be useful to show the properties of the numerical kinetic scheme such as the
preservation of the steady state water at rest, and the positivity of the wetted area.

Section \ref{kineticscheme} is devoted to the construction of the kinetic scheme. The upwin\-ding of
the source term due to the topography
is performed in a close manner described by  Perthame et al. \cite{PS01} using an energetic balance
at microscopic level
for the Shallow Water equations.

Finally,  we present in section \ref{numerics} a numerical validation of this study by the
comparison between the resolution of this model and the resolution of the Allievi equation solved by
the research code \verb+belier+ used at Center in Hydraulics Engineering of Electricit\'{e} De
France (EDF) \cite{W93} for the case of critical waterhammer tests.

%%%%%%%%%%%%%%%%%%%%%%%%%%%%%%%%%%%%%%%%%%%%%%%%%%%%%%%%%%%%%%%%
\section{The mathematical model and the kinetic formulation}\label{model}
%%%%%%%%%%%%%%%%%%%%%%%%%%%%%%%%%%%%%%%%%%%%%%%%%%%%%%%%%%%%%%%%

We derived a conservative model for pressurised
flows from the 3D system of compressible Euler equations by integration
over sections orthogonal to the flow axis.
\subsection{The mathematical model : a ``Shallow Water like'' system of equations}
The equation for conservation of mass and the first equation
for the conservation of momentum are:
\begin{eqnarray}
\partial_t \rho + \div(\rho \,\vec U)	 & = & 0 \label{E1}\\
\partial_t (\rho \,u) + \div(\rho \,u\,\vec U) & =  & F_x  -\partial_x P \label{E2}
\end{eqnarray}
with the speed vector $\vec U = u\vec i + v\vec j + w\vec k =u\vec i + \vec V$, where the unit
vector $\vec i$ is
along the main axis, $\rho$ is the density of the water.
We use	the Boussinesq linearised pressure law (see \cite{SWB98}):
\begin{equation*}
P = P_a  + \frac{1}{\beta }\left( {\frac{\rho}{{\rho _0 }} - 1} \right),
\end{equation*}
where $\rho _0$ is the density at the atmospheric pressure $P_a$ and $\beta$ the coefficient of
compressibility of the water.
Exterior strengths $\vec F$ are the gravity $\vec g$ and
the friction term $S_f$ which is assumed to be given by the Manning-Strickler law (see
\cite{SWB98}):
\begin{equation} \label{strickler}
S_f=K \,u\,\abs{u}\quad\hbox{ with }\quad K = \frac{1}{K_s^2\,R_h^{4/3}}
\end{equation}
where $K_s>0$ is the Strickler coefficient, depending on the material,	and $R_h$ is the so called
hydraulic
radius given by $R_h=\displaystyle \frac{S}{P_m}$.
$S$ represents the cross-section area of the pipe whereas $P_m$ is the perimeter of the section.
Then Equations (\ref{E1})-(\ref{E2}) become:
\begin{eqnarray*}
\partial _t \rho  + \partial _x (\rho\, u) + \div_{(y,z)} (\rho\, \vec V)&=&0 \label{E'1}\\
\partial _t (\rho\, u) + \partial _x (\rho\, u^2 ) + \div_{(y,z)} (\rho\, u\,\vec V)& = &-\rho
g(\partial_x Z + S_f )
- \frac{\partial_x\rho}{\beta\rho_0} \quad .\label{E'2}
\end{eqnarray*}
Assuming that the pipe is infinitely rigid and has a uniform constant cross-section $S$, and taking
averaged values in sections
orthogonal to the main flow axis, we get the following system written in a conservative form for the
unknowns
$M = \rho \,S \,,\, D = \rho\, S \, u$:
 \begin{eqnarray}
\label{masseCH} \partial_t (M)+\partial_x (D) & = & 0\\
\label{qmouvCH} \partial_t (D)+\partial_x\left(\frac{D^2}{M}+c^2 \,M \right)& = & - g \,M (\partial_x Z + S_f)
\end{eqnarray}
where $\displaystyle c=\frac{1}{\sqrt{\beta\,\rho_0}}$ is the speed of sound.
A complete derivation of this model, taking into account the deformations of the pipe, contracting
or expanding sections, and a spatial second order Roe-like finite volume method in a  linearly
implicit version is presented in \cite{BG08} (see \cite{BG02} for the  first order implicit scheme).
This system of partial differential equation is formally close to the Shallow Water equations where
the conservative variables are the wet area and the discharge, thus we define an ``FS-equivalent''
wet area (FS for Free Surface) $A$  and a ``FS-equivalent discharge'' $Q$ through the relations:
$$M=\rho\,S=\rho_0\,A \quad \mbox{and} \quad D=\rho\, S \, u =\rho_0\,Q \quad .$$
Dividing (\ref{masseCH})-(\ref{qmouvCH}) by $\rho_0$ we can write this system
under the conservative form:
\begin{equation}\label{CHconservative}
\partial_t U+\partial_x F(U)=G(x,U)
\end{equation}
where the unknown state is
$U=(A,Q)^t$, the flux vector is
$F(U)=(Q,\displaystyle \frac{Q^2}{A}+ c^2 A)^t$
and the source term writes
$G(x,U)=(0,\displaystyle -g A(\partial_x Z + S_f))^t$. This new set of variables allows a more
natural treatment of mixed flows (see \cite{BG08}).

Let us now recall the main properties of the system (\ref{CHconservative}) whose proofs
can be found in
\cite{BGG08}.
\begin{thm}\label{CHhyp}
The system (\ref{CHconservative}) is strictly hyperbolic. It admits a mathematical
entropy:
\begin{equation}\label{entroCH}
E(A,Q,Z) = \frac{Q^2}{2 A} + g A Z + c^2 A \ln A
\end{equation}
which satisfies, for smooth enough solution, the entropy inequality :
\begin{equation*}
\partial_t E + \partial_x [u (E + c^2 A)] \leq 0 \quad .
\end{equation*}
Also, for the frictionless pipes ($S_f = 0$), the system (\ref{CHconservative}) admits a family of smooth 
steady states characterized by the relations:
\begin{equation*}
Q = A u = C_1 \;,
\end{equation*}
\begin{equation*}
\frac{u^2}{2}  +  g \, Z + c^2 \ln A= C_2 \;,
\end{equation*}
where $C_1$ and $C_2$ are two arbitrary constants. The quantity $\displaystyle \frac{u^2}{2} + g \,
Z + c^2 \ln A$
is also called the total head.
\end{thm}
\begin{rem} \rm 
An easy computation leads to the equality:
\begin{equation*}
\partial_t E + \partial_x [u (E + c^2 A)]  = - u S_f = -K \vert u \vert^3 \leq 0 \; .
\end{equation*}
Thus for a frictionless pipe we obtain an entropy equality whereas the entropy inequality is strict as soon as
a friction term is considered.
 
Let us also remark that the still water steady state for frictionless pipe, namely $u \equiv 0$,  satisfies:
$ g \, Z + c^2 \ln A = C_2$.
\end{rem}
%%%%%%%%%%%%%%%%%%%%%%%%%%%%%%%%%%%%%%%%%%%%%%%%%%%%%%%%%%%%%%%%%%
\subsection{The kinetic approach}
%%%%%%%%%%%%%%%%%%%%%%%%%%%%%%%%%%%%%%%%%%%%%%%%%%%%%%%%%%%%%%%%%%
We present in this section the kinetic formulation for pressurised flows in closed water pipes
modelised by the preceding system
of partial differential equations (see \cite{BGG08} for more details and properties).
Let us mention that the following results (namely Theroem \ref{CHcin} and Theorem \ref{energyCH}) are only valid for 
frictionless pipes.
Let us consider a smooth real function $\chi$ which has the following properties:
\begin{equation}\label{propchi}
\chi(\omega)=\chi(-\omega) \geq 0\;,\;
\int_\R \chi(\omega) d\omega =1,
\int_\R \omega^2 \chi(\omega) d\omega=1 \;.
\end{equation}
We then define the density of particles  $\mathcal {M}(t,x,\xi)$ by the so-called {\it Gibbs
equilibrium}:
$$\mathcal {M}(t,x,\xi) =\frac{A(t,x)}{c} \chi \left( \frac{\xi-u(t,x)}{c} \right)
\quad .$$
These definitions allow to obtain a kinetic representation of the system
(\ref{CHconservative})
by the following result (see \cite{BGG08} for the proof).
%%%%%%%%%%%%%%%%%%%%%%%%%%%%%%%%%%%%%%%%%%%%%%%%%%%%%%%%%%%%%%%%%%
% Formulation cinetique pour CH
%%%%%%%%%%%%%%%%%%%%%%%%%%%%%%%%%%%%%%%%%%%%%%%%%%%%%%%%%%%%%%%%%%
\begin{thm}\label{CHcin}
The couple of functions $(A,Q)$ is a strong solution of the system
(\ref{CHconservative})	if and only if
${\mathcal M}$ satisfies the kinetic equation:
\begin{equation}\label{CHeqcin}
\frac{\partial}{\partial t} {\mathcal M} +\xi \cdot \frac{\partial}{\partial x} {\mathcal M}
-g\frac{\partial}{\partial x}Z \cdot \frac{\partial}{\partial
\xi} {\mathcal M} =K(t,x,\xi)
\end{equation}
for some collision term $K(t,x,\xi)$ which satisfies for a.e. $(t,x)$
\begin{equation*}\label{collisionch}
\displaystyle \int_\R K\, d\xi=0 \;,\; \displaystyle \int_\R  \xi\, K d\, \xi=0 \quad .
\end{equation*}
\end{thm}
This result is a consequence of the following relations verified by the microscopic equilibrium:
\begin{eqnarray}
\label{macrochA} A &=& \int_{\R} {\mathcal M}(\xi) \,d\xi \;,\\
\label{macrochQ} Q &=& \int_{\R} \xi {\mathcal M}(\xi)\, d\xi \;,\\
\label{macrochI} \frac{Q^2}{A} + c^2 A&=& \int_{\R} \xi^2 {\mathcal M}(\xi) \,d\xi \;.
\end{eqnarray}

This theorem produces a very useful consequence: the nonlinear system
(\ref{CHconservative}) can be viewed as
a simple linear equation on a nonlinear quantity $\mathcal{M}$ for which it is easier to find simple
numerical schemes with good theoretical properties: it is this feature which will be exploited to
construct a kinetic scheme.

\begin{thm}\label{energyCH}
Let $A(x,t) > 0$ and $Q(x,t)$ be two given functions.
\begin{enumerate}
\item The minimum of the energy:
$$\mathcal {E}(f)=\int_{\R} \left(
\frac{\xi^2}{2}f(\xi)+c^2 f(\xi) \ln(f(\xi)) +g Z f(\xi) +c^2 \ln(c\sqrt{2
\pi}) f( \xi)\right)d\xi \quad ,$$
under the constraints:
 $$f \geq 0\,, \; \int _\R f(\xi) d \xi= A\,, \ \int_\R \xi f(\xi) d\xi=Q \quad ,$$
is attained by the function:
$$\displaystyle \mathcal {M}(t,x,\xi) = \frac{A}{c} \chi \left(\frac{\xi-u(t,x)}{c} \right)$$
% \mathcal {M}(A,\xi-u)=
where $\chi$ is defined by:
\begin{equation}\label{choixCH}
\chi(\omega)=\displaystyle \frac{1}{\sqrt{2 \pi}} \exp\left(-\frac{\omega^2}{2}\right) \quad .
\end{equation}
\item Moreover, the function $\chi$ defined by (\ref{choixCH}) ensures us to have the relation
$$\mathcal{E}(\mathcal{M}) = E(A,Q,Z)$$
if $A$ and $Q$ are solution of the pressurised flow equations
(\ref{CHconservative})
and the entropy $E$ is defined by (\ref{entroCH}).
\item The Gibbs equilibrium $\mathcal{M}$ satisfies the still water steady state
equation.
\end{enumerate}
\end{thm}
%%%%%%%%%%%%%%%%%%%%%%%%%%%%%%%%%%%%%%%%%%%%%%%%%%%%%%%%%%%%%%%%%%%%%%%%%%%%%%%%%%%%
\section{The kinetic scheme}\label{kineticscheme}
%%%%%%%%%%%%%%%%%%%%%%%%%%%%%%%%%%%%%%%%%%%%%%%%%%%%%%%%%%%%%%%%%%%%%%%%%%%%%%%%%%%%
The spatial domain is a pipe of length $L$. The main axis of the pipe is divided in $N$ meshes
$\displaystyle m_i=]x_{i-1/2},x_{i+1/2}[, \ 1\leq i\leq N \,,\, \mbox{ of length } h_i \mbox{ and
center } x_i$. We denote
$\Delta x = \min_{1\leq i\leq N} h_i$. $\Delta t$ denotes the time step at time $t_n$ and we set
$t_{n+1}=t_n + \Delta t$.\\
The discrete macroscopic unknowns are $U_i^n=\left(\begin{array}{c}A_i^n\\
Q_i^n\end{array}\right)$ with $1\leq i\leq N$ and $ 0\leq n \leq n_{max}$. They represent the mean value
of
$U$ on the cell $m_i$ at time $t_n$.

If $Z(x)$ is the function describing the bottom elevation, its piecewise constant representation is
given by
$\bar{Z}(x) = Z_i \1_{m_{i}}(x)$ with $Z_i = Z(x_i)$ for example.

Replacing $Z$ by $\bar{Z}$ and neglecting the collision term $K(t,x,\xi)$ in a first step, the
Equation (\ref{CHeqcin}) in the cell $m_i$ writes:
%  Equation in the cell m_i
\begin{equation}\label{eqcin}
\frac{\partial}{\partial t} {\mathcal M} +\xi \cdot \frac{\partial}{\partial x} {\mathcal M}
=0 \quad \mbox{ for } x \in m_i \quad .
\end{equation}
This equation is a linear transport equation whose explicit discretisation may be done directly by
the following way.
Denoting for  $x\in m_i \,,\, f(t_n,x,\xi)~=~\mathcal{M}_i^n(\xi) $\linebreak%={M}(A_i^n,Q_i^n,\xi)$
the maxwellian state associated to  $A_i^n \,,\mbox{ and } Q_i^n$, the usual finite volume discretisation of the Equation (\ref{eqcin}) leads to:
%  Equation in the cell m_i : discretisation
\begin{equation} \label{cindiscret}
f_i^{n+1}(\xi) = \mathcal{M}_i^n(\xi)+\frac{\Delta t}{h_i}\,
\xi \, \left(\mathcal{M}_{i+\frac{1}{2}}^-(\xi)-\mathcal{M}_{i-\frac{1}{2}}^+(\xi)\right)
\end{equation}
where the fluxes $\mathcal{M}_{i+\frac{1}{2}}^\pm$ have to take into account the discontinuity of
the altitude $\bar{Z}$ at the cell interface $x_{i+1/2}$.
Indeed, noticing that the fluxes can also be written as:
% Remark on the difference on the fluxes
\begin{equation*}
\mathcal{M}_{i+\frac{1}{2}}^-(\xi) = \mathcal{M}_{i+\frac{1}{2}} +
\left(\mathcal{M}_{i+\frac{1}{2}}^- - \mathcal{M}_{i+\frac{1}{2}}\right)
\end{equation*}
the quantity $\delta \mathcal{M}_{i+\frac{1}{2}}^- = \mathcal{M}_{i+\frac{1}{2}}^- -
\mathcal{M}_{i+\frac{1}{2}}$
holds for the discrete contribution of the source term $g A \partial_x Z$ in the system for negative
velocities $\xi \leq 0$  due to the upwinding of the source term.
Thus $\delta \mathcal{M}_{i+\frac{1}{2}}^-$ has to vanish for positive velocity $\xi > 0$,
as proposed by the choice of the interface fluxes below.
% The choice for the fluxes.

\begin{figure}[htp]
\centering{
\includegraphics[scale=0.5]{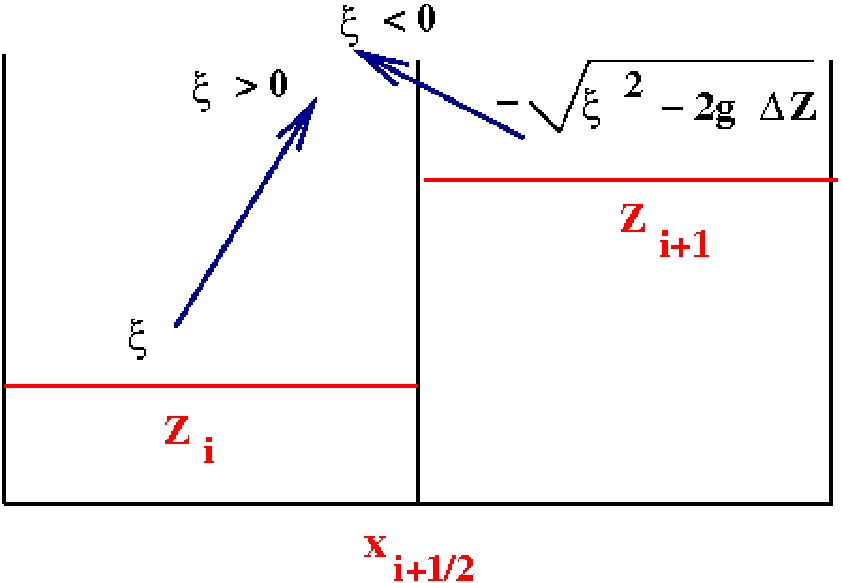}
}
\caption{Transmission}\label{transmission}
\end{figure}

\begin{figure}[htp]
\centering{
\includegraphics[scale=0.5]{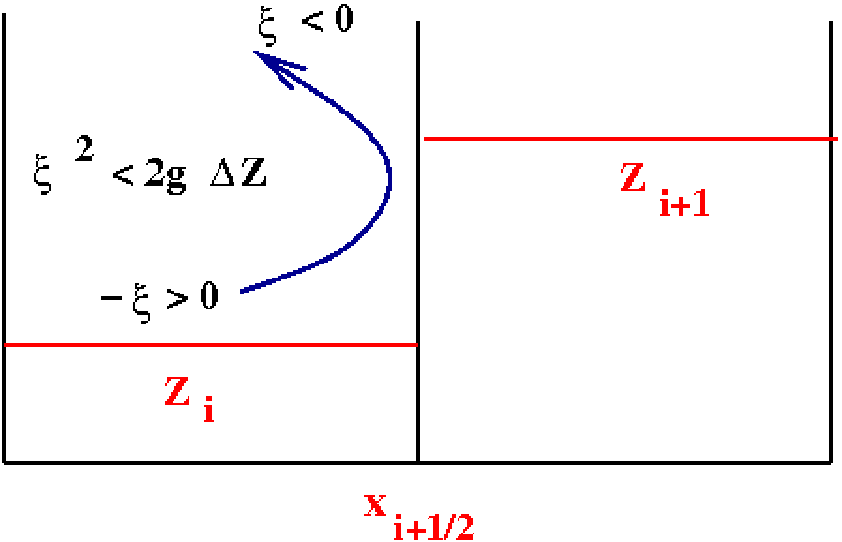}
}
\caption{Reflection}\label{reflection}
\end{figure}

Let us now detail our choice for the fluxes $\mathcal{M}_{i+\frac{1}{2}}^\pm$ at the interface. It
can be justified by using a generalised characteristic method for the Equation (\ref{CHeqcin})
(without the collision kernel) but we give instead a presentation based on some physical energetic
balance.
Let us denote  $\displaystyle \Delta^-Z_{i+\frac{1}{2}}=Z_{i+1}-Z_i$ and
$\displaystyle\Delta^+Z_{i+\frac{1}{2}}=Z_i-Z_{i+1}$.
In order to take into account the neighboring cells by means of a natural interpretation of the
microscopic features of the system, we formulate a peculiar  discretisation for the fluxes in
(\ref{cindiscret}), computed by the following upwinded formulas:
% The formulas !!!
\begin{eqnarray} \label{Mmoins}
\mathcal{M}^-_{i+\frac{1}{2}} (\xi) & =&	\mathcal{M}_i^n(\xi)\,\1_{\xi\geq 0}
+ \overset{reflection} {\overbrace{ \mathcal{M}_i^n(-\xi)\,\1_{\xi^2\leq
2g\Delta^-Z_{i+\frac{1}{2}}}\,\1_{\xi\leq 0} }}\\
& &+\underset{transmission} {\underbrace{
\mathcal{M}_{i+1}^n\left(-\sqrt{\xi^2-2g\Delta^-Z_{i+\frac{1}{2}}}\right)
\,\1_{\xi^2\geq2g\Delta^-Z_{i+\frac{1}{2}}}\,\1_{\xi\leq 0} }}\nonumber
\end{eqnarray}
\begin{eqnarray}\label{Mplus}
\mathcal{M}^+_{i+\frac{1}{2}}(\xi) &=& \mathcal{M}_{i+1}^n(\xi)\,\1_{\xi\leq 0}
+ \overset{reflection} {\overbrace{ \mathcal{M}_{i+1}^n(-\xi)\,\1_{\xi^2\leq
2g\Delta^+Z_{i+\frac{1}{2}}}\,\1_{\xi\geq 0} }}\\
&&+\underset{transmission} {\underbrace{
\mathcal{M}_i^n\left(\sqrt{\xi^2-2g\Delta^+Z_{i+\frac{1}{2}}}\right)
\,\1_{\xi^2\geq 2g\Delta^+Z_{i+\frac{1}{2}}}\,\1_{\xi\geq 0} }}\nonumber
\end{eqnarray}
% Explanation
The effect of the source term is made explicit by treating it as a physical potential. The choices
(\ref{Mmoins})-(\ref{Mplus})
are thus a mathematical formalization to describe the physical microscopic behaviour of the system.
The contribution of the interface $x_{i+1/2}$ to
$f_i^{n+1}$ is	given by:
\begin{itemize}
 \item the particles in the cell $m_i$ at time $t_n$ with non negative velocities $\xi$ through the
term $\mathcal{M}_i^n(\xi)\,\1_{\xi\geq 0}$ and those of them that are reflected  (thus taken into account with velocity
$-\xi$) if their kinetic energy is not large enough to overpass
the potential  difference i.e. $\xi^2\leq 2g\Delta^\pm
Z_{i+\frac{1}{2}}$: see Figure \ref{reflection}.
\item the particles in the cell $m_{i+1}$  at time $t_n$ with a kinetic
energy enough to overpass the potential  difference ( $\xi^2\geq 2g\Delta^\pm Z_{i+\frac{1}{2}}$ )
and speed up or down  according to this potential jump.
It is the transmission phenomenon in classical mechanics as shown in Figure
\ref{transmission}.
\end{itemize}
Since we neglected the collision term, it is clear that $f^{n+1}$ computed by the discretised
kinetic equation (\ref{cindiscret})
is no more a Gibbs equilibrium. Therefore, to recover the macroscopic variables $A$ and $Q$, according to  the identities
(\ref{macrochA})-(\ref{macrochQ}), we set:
\begin{equation}\label{macro}
 U_i^{n+1}=
\left(\begin{array}{l}
A_i^{n+1} \\
Q_i^{n+1}
\end{array}
\right)
\overset{def}{=} \int_\R \left(\begin{array}{l}
1 \\
\xi
\end{array}
\right) f_i^{n+1} \,d\xi
\end{equation}
Now, we can integrate the discretised kinetic equation (\ref{cindiscret}) against 1 and $\xi$ to
obtain the macroscopic kinetic scheme:
\begin{equation} \label{macrodiscret}
U_i^{n+1} = U_i^n +\frac{\Delta t}{h_i}\left(F_{i+\frac{1}{2}}^- - F_{i-\frac{1}{2}}^+ \right)
\end{equation}
The numerical fluxes are thus defined by the kinetic fluxes as follows:
\begin{equation}\label{fluxdiscret}
F_{i+\frac{1}{2}}^\pm
\overset{def}{=} \int_\R \xi \left(\begin{array}{l}
1 \\
\xi
\end{array}
\right) \mathcal{M}^\pm_{i+\frac{1}{2}}(\xi)\,d\xi
\end{equation}
\begin{rem} \label{remnum} \rm ~
\begin{itemize}
\item We see immediately that the kinetic scheme (\ref{macrodiscret})-(\ref{fluxdiscret}) is wetted area conservative.
Indeed, let us denote the first component of the discrete fluxes (\ref{fluxdiscret})
$\left(F_A\right)^\pm_{i+\frac{1}{2}}$:
\begin{equation*}
\left(F_A\right)^\pm_{i+\frac{1}{2}}
\overset{def}{=} \int_\R \xi \mathcal{M}^\pm_{i+\frac{1}{2}}(\xi)\,d\xi
\end{equation*}
An easy computation using the change of variables $\mu = \abs{\xi}^2 - 2g\Delta^+ Z_{i+\frac{1}{2}}$
in the formulas
(\ref{Mmoins})-(\ref{Mplus}) defining the kinetic fluxes $\mathcal{M}_{i+\frac{1}{2}}^\pm$ allows us
to show that:
\begin{equation*}
\left(F_A\right)^+_{i+\frac{1}{2}} = \left(F_A\right)^-_{i+\frac{1}{2}}
\end{equation*}
\item Computing the macroscopic state $U$ by the formula (\ref{macro}) or the fluxes by the formula
(\ref{fluxdiscret}) is not easy if the function $\chi$ verifying
the properties (\ref{propchi}) is not compactly supported.
We use instead the function defined by:
\begin{equation} \label{chichoix}
\chi(\omega)=\frac{1}{2\,\sqrt{3}}\,\1_{[-\sqrt{3},\sqrt{3}]}(\omega) \; .
\end{equation}
We get 
$\mathcal{M}_i^n(\xi)=\displaystyle\frac{A_i^n}{2\,c\,\sqrt{3}}\,\1_{[u_i^n-c\,\sqrt{3},u_i^n+c\,\sqrt{3}]}(\xi) .$
Of course, from Theorem \ref{energyCH}, the property on the microscopic energy and the still water steady state
is no more valid but we will prove in Theorem \ref{prop} that our proposed kinetic scheme preserves the still water
steady state (for the frictionless pipes) and the positivity of the equivalent wetted area.
\item In the case where the friction $S_f$ defined by (\ref{strickler}) is present, it is added at the macroscopic level
(\ref{macrodiscret}) as an extra source term.
\end{itemize}
\end{rem}
We are now able to state the main properties of the kinetic scheme.
\begin{thm}\label{prop}
We choose the function $\chi$ defined by the formula (\ref{chichoix}) and we assume the CFL condition
\begin{equation} \label{CFL}
\Delta t \max_{1\leq i\leq N}\left(\abs{u^n_i} +c \sqrt{3} \right) \leq \Delta x.
\end{equation}
Then
\renewcommand{\labelenumi}{(\roman{enumi})}
\begin{enumerate}
\item the kinetic scheme (\ref{macrodiscret})-(\ref{fluxdiscret}) keeps the pseudo wetted area $A^n_i$
positive.
\item the kinetic scheme  (\ref{macrodiscret})-(\ref{fluxdiscret}) preserves the still water steady
state,
$$ u^n_i = 0 \,,\, g \, Z_i + c^2 \ln A_i = K$$
\end{enumerate}
\end{thm}
{\bf Proof of theorem \ref{prop}}
Let us mention that the  CFL condition (\ref{CFL}) is obtained from the linear discretised kinetic transport
equation  (\ref{cindiscret}) for the particular choice of the function  $\chi$ defined by the formula (\ref{chichoix}).
This condition ensures the positivity of the pseudo wetted area as we show below.

Since $A_i =\displaystyle \int_\R
f_i^{n+1} \,d\xi$, it is sufficient to prove that
$f_i^{n+1} \geq 0$. Writing the microscopic scheme  (\ref{cindiscret}), (\ref{Mmoins}), (\ref{Mplus}),
using the CFL condition (\ref{CFL}), and the fact that the function $\chi$ that we have chosen is compactly supported,
one may see that if we suppose that
$A_i^n \geq 0$, then $f_i^{n+1}$ is a sum of non-negative quantities. For the second point, setting
$u_i^n = 0$, we prove easily
that in the discretised kinetic equation (\ref{cindiscret}), we have
$$\mathcal{M}^-_{i+\frac{1}{2}} (\xi)  = \mathcal{M}^+_{i-\frac{1}{2}} (\xi) .$$
This implies $f^{n+1}_i = \mathcal{M}^n_{i} (\xi)$, which ensures by definition
$A^{n+1}_i= A^{n}_i$ and $Q^{n+1}_i=Q^{n}_i$. Thus we obtain $u^{n+1}_i = 0$.

%%%%%%%%%%%%%%%%%%%%%%%%%%%%%%%%%%%%%%%%%%%%%%%%%%%%%%%%%%%%%%%%%%%%%%%%%%%%%%%%%%%%%%%%%%%%%%%%%%%%
\section{Numerical validation}\label{numerics}
%%%%%%%%%%%%%%%%%%%%%%%%%%%%%%%%%%%%%%%%%%%%%%%%%%%%%%%%%%%%%%%%%%%%%%%%%%%%%%%%%%%%%%%%%%%%%%%%%%%%

We present now numerical results of a water hammer test.
The pipe of circular cross-section of $2 \mbox{ m}^2$ (the diameter therefater is denoted by $\delta$)
and thickness $e=20$ cm is $2000$ m long. The altitude of the upstream end of the
pipe is $250$ m and the slope is $5^{\circ}$. 
The Young modulus is $23 \, 10^9 \mbox{ Pa}$ since the
pipe is supposed to be built in concrete. \\
The density at the atmospheric pressure $\rho _0$ is  $1000 \mbox{ kg}/\mbox{m}^3$ and the coefficient of
compressibility of the water $\beta$ is $5.0 \, 10^{-10} \mbox{ Pa}^{-1}$. \\
The wave speed is thus obtained by the formula (see \cite[formula (2.39)]{SWB98}):
\begin{equation}
 a =\displaystyle \frac{c}{\sqrt{1+\displaystyle\frac{\delta}{\beta\,e\,E}}} = 1086.6 \mbox{ m}/\mbox{s}^{-1} \quad .
\end{equation}
The total upstream head is 300 m. The initial downstream discharge is $10 \mbox{ m}^3/\mbox{s}$
and we cut the flow in $5$ seconds. Let us define the
piezometric line by:
\begin{equation} \label{piezo}
\displaystyle piezo = z + \delta + p \; \mbox{ with } p = \frac{c^2 \, (\rho -\rho_0)}{\rho_0 \,
g} \: .
\end{equation}
We present a validation of the proposed scheme by comparing numerical results of the proposed model
solved by the kinetic scheme with the ones obtained by solving Allievi equations by the method of
characteristics with the so-called \verb+belier+ code: an industrial code used by the engineers of the
Center in Hydraulics Engineering of Electricit\'{e} De France (EDF) \cite{W93}.
Our code is written in Fortran and runs during a few seconds on LinuX, Windows and MacIntosh operating systems.

\begin{figure}[htp]
\centering{
\includegraphics[scale=0.95]{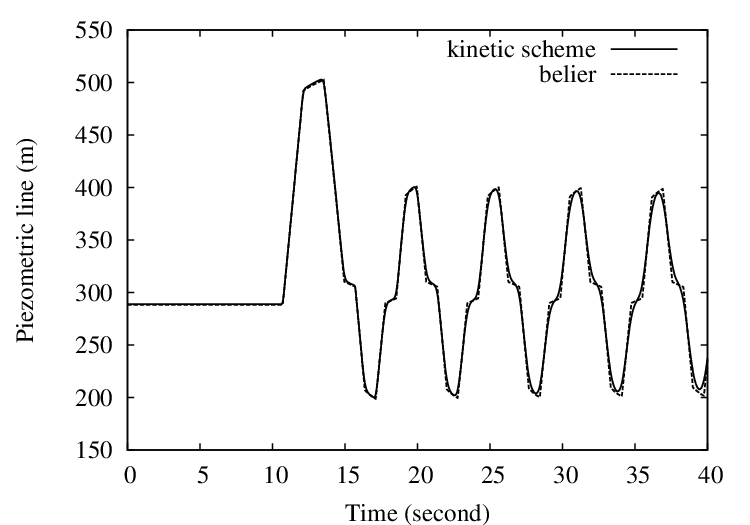}
\includegraphics[scale=0.95]{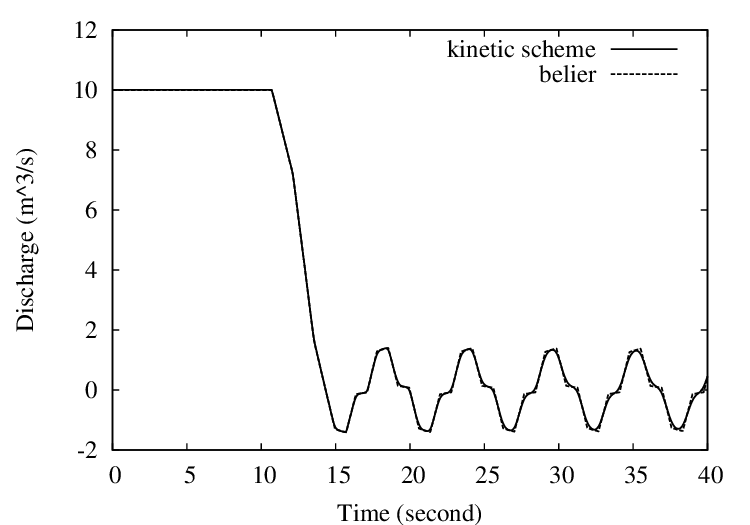}
\caption{Comparison between the kinetic scheme and the industrial code belier}
\label{order1}
{\small Piezometric line (top) and
discharge  (bottom) at the middle of the pipe}
}
\end{figure}

\begin{figure}[htp]
\centering{
\includegraphics[scale=.95]{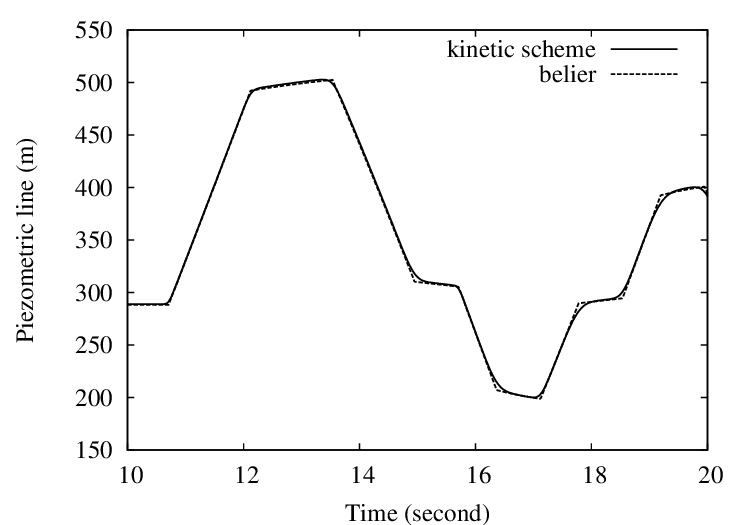}
\includegraphics[scale=.95]{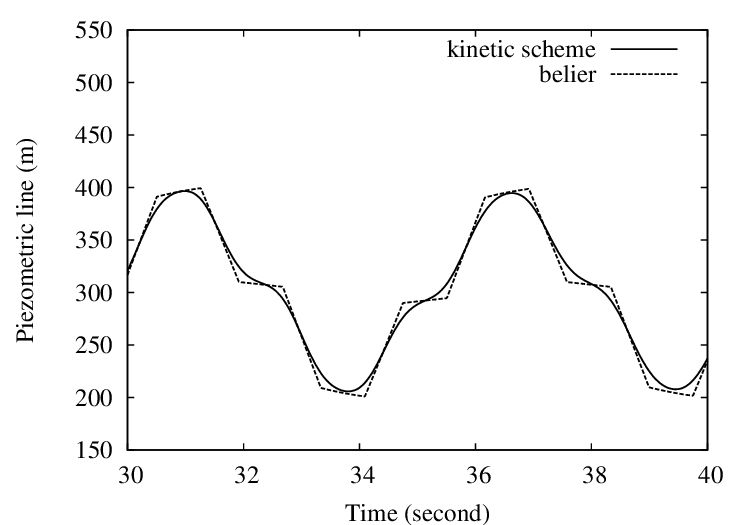}
\caption{Comparison between the kinetic scheme and the industrial code belier}
\label{zoom}
{\small Beginning of the simulation (top) and
end of the simulation (bottom) at the middle of the pipe}
}
\end{figure}

A simulation of the water hammer test was done for a CFL coefficient equal to 0.8 (i.e. $CFL = 0.8$)
and a spatial discretisation of 1000 mesh points (the mesh size is equal to $2$ m).
In the Figure \ref{order1}, we present a comparison between the results obtained by our kinetic scheme and
the ones obtained by the ``belier'' code at the middle of the pipe:
the behavior of the piezometric line, defined by Equation (\ref{piezo}), and the discharge at the middle of the pipe.
One can observe that the results for the proposed model are in very good agreement with the  solution of
Allievi equations. In Figure \ref{zoom}, we present the piezometric line for the beginning and the end of the simulation.
One can see that the peak of pressure (observed in the beginning of the simulation) is very well obtained by the kinetic scheme. 
Thus the strength of the water hammer is very good predicted by the proposed numerical scheme. 
A little smoothing effect, observed at the end of the simulation,  may be probably due to the first order discretisation type.
A second order scheme could be implemented naturally and produce a better approximation.
\section{Conclusion}
As mentionned in the introduction, our goal is to build a kinetic scheme for mixed flows.
Perthame et al. \cite{BPV03,PS01} have shown that the kinetic approach is relevant for free surface flows in open
channels: the resulting kinetic scheme is easily implemented an enjoyed very good properties (posivity of the wetted area 
and discrete entropy inequalities). For pressurised flows, we have shown that this approach is also relevant.
This allows us to investigate the construction of a kinetic scheme for mixed flows in closed water pipes.
\bibliographystyle{plain}
%\bibliography{cinemix}

\end{document}